\documentclass[conference]{IEEEtran}
\usepackage{times,amsmath,epsfig}
\usepackage[latin5]{inputenc}
\DeclareGraphicsExtensions{.wmf,.eps}

\begin{document}

\title{Linear Differential Equations\\ with Fuzzy Boundary Values}

\author{\IEEEauthorblockN{Nizami Gasilov}
\IEEEauthorblockA{Baskent University,\\Eskisehir yolu 20. km, Baglica,\\06810 Ankara, Turkey\\Email: gasilov@baskent.edu.tr}
\and \IEEEauthorblockN{Şahin Emrah Amrahov}
\IEEEauthorblockA{Ankara University,\\Computer Engineering
Department,\\06100 Ankara, Turkey\\Email:
emrah@eng.ankara.edu.tr} \and \IEEEauthorblockN{Afet Golayoğlu
Fatullayev} \IEEEauthorblockA{Baskent University,\\Eskisehir
yolu 20. km, Baglica,\\06810 Ankara, Turkey\\Email:
afet@baskent.edu.tr}}
 \maketitle

\begin{abstract}
In this study, we consider a linear differential equation with fuzzy boundary values. We express the solution of the problem in terms of a fuzzy
set of crisp real functions. Each real function from the solution set 
satisfies differential equation, and its boundary values belong to
intervals, determined by the corresponding fuzzy numbers. The least
possibility among possibilities of boundary values in corresponding fuzzy
sets is defined as the possibility of the real function in the fuzzy
solution.

In order to find the fuzzy solution we propose a method based on the properties of linear transformations. We show that, if the corresponding crisp problem has a unique solution then the fuzzy problem has unique
solution too. We also prove that if
the boundary values are triangular fuzzy numbers, then the value of the
solution at any time is also a triangular fuzzy number.

We find that the fuzzy solution determined by our method is the same as the one that is obtained from solution of crisp problem by the application of the extension principle.

We present two examples describing the proposed
method.
\end{abstract}



\IEEEpeerreviewmaketitle \noindent \bigskip

\section{Introduction}

Approaches to fuzzy boundary value problems can be of two types. The first
approach assumes that, even if only the boundary values are fuzzy in the
handling problem, the solution is fuzzy function, consequently, the
derivative in the differential equation can be considered as a derivative of
fuzzy function. This derivative can be Hukuhara derivative, or a derivative
in generalized sense. Bede \cite{Bede06} has demonstrated that a large class
of boundary value problems have not a solution, if Hukuhara derivative is
used. To overcome this difficulty, in \cite{BG05} and \cite{CCRF06} the
concept of generalized derivative is developed and fuzzy differential
equations have been investigated using this concept (see also \cite%
{BRB07,CCRFRM08,CCRMRF07}). Recently, Khastan and Nieto \cite{KhN10} have
found solutions for a large enough class of boundary value problems with the
generalized derivative. However as it is seen from the examples in mentioned
article, these solutions are
difficult to interpret because four different problems, obtained by using
the generalized second derivatives, often does not reflect the nature of the
problem.

The second approach is based on generating the fuzzy solution from the crisp
solution. In particular case, for the fuzzy initial value problem this approach
can be of three ways. The first one uses the extension principle. In this
way, the initial value is taken as a real constant, and the resulting crisp
problem is solved. Then the real constant in the solution is replaced with the initial fuzzy value. In the final solution, arithmetic
operations are considered to be operations on fuzzy numbers (\cite{BF00,BF01}%
). The second way, offered by Hüllemerier \cite{Hul97}, uses the concept of
differential inclusion. In this way, by taking an alpha-cut of initial
value, the given differential equation is converted to a differential
inclusion and the obtained solution is accepted as the alpha-cut of the
fuzzy solution. Misukoshi et al \cite{MBCRB07} have proved that, under
certain conditions, the two main ways of the approach are equivalent for the
initial value problem. The third way is offered by Gasilov et al \cite{GAF11}%
. In this way the fuzzy problem is considered to be a set of crisp problems.

In this study, we investigate a differential equation with fuzzy boundary
values. We interpret the problem as a set of crisp problems. For linear
equations, we propose a method based on the properties of linear
transformations. We show that, if the solution of the corresponding crisp
problem exists and is unique, the fuzzy problem also has unique solution.
Moreover, we prove that if the boundary values are triangular fuzzy numbers,
then the value of the solution is a triangular fuzzy number at each time. We
explain the proposed method on examples.

\section{Fuzzy Boundary Value Problem}

Below, we use the notation $\widetilde{u}=(u_{L}(r),u_{R}(r))$, $(0\leq
r\leq 1)$, to indicate a fuzzy number in parametric form. We denote $%
\underline{u}=u_{L}(0)$ and $\overline{u}=u_{R}(0)$ to indicate the left and
the right limits of $\widetilde{u}$, respectively. We represent a triangular
fuzzy number as $\widetilde{u}=(l,m,r)$, for which we have $\underline{u}=l$%
\ and $\overline{u}=r$.

In this paper we consider a fuzzy boundary value problem (FBVP) with crisp
linear differential equation but with fuzzy boundary values. For clarity we
consider second order differential equation: 
\begin{equation}
\left\{ 
\begin{array}{c}
x^{\prime \prime }+a_{1}(t)x^{\prime }+a_{2}(t)x=f(t) \\ 
x(0)=\widetilde{A} \\ 
x(T)=\widetilde{B}%
\end{array}%
\right.  \label{G_FBVP}
\end{equation}

\noindent We note that the coefficients of the differential equation are not
necessary constant.

Let us represent the boundary values as $\widetilde{A}=a_{cr}+\widetilde{a}$%
\ and $\widetilde{B}=b_{cr}+\widetilde{b}$, where $a_{cr}$ and $b_{cr}$ are
vectors with possibility 1 and denote the crisp parts (the vertices) of $%
\widetilde{A}$ and $\widetilde{B}$; $\widetilde{a}$ and $\widetilde{b}$
denote the uncertain parts with vertices at the origin.

We split the problem (\ref{G_FBVP}) to following two problems:

1) Associated crisp problem (which is non-homogeneous): 
\begin{equation}
\left\{ 
\begin{array}{c}
x^{\prime \prime }+a_{1}(t)x^{\prime }+a_{2}(t)x=f(t) \\ 
x(0)=a_{cr} \\ 
x(T)=b_{cr}%
\end{array}%
\right.  \label{N-H_CBVP}
\end{equation}

2) Homogeneous problem with fuzzy boundary values: 
\begin{equation}
\left\{ 
\begin{array}{c}
x^{\prime \prime }+a_{1}(t)x^{\prime }+a_{2}(t)x=0 \\ 
x(0)=\widetilde{a} \\ 
x(T)=\widetilde{b}%
\end{array}%
\right.  \label{FBVP}
\end{equation}

It is easy to see that, a solution of the given problem (\ref{G_FBVP}) is of
the form $x(t)=x_{cr}(t)+x_{un}(t)$ (crisp solution + uncertainty). Here $%
x_{cr}(t)$ is a solution of the non-homogeneous crisp problem (\ref{N-H_CBVP}%
); while $x_{un}(t)$ is a solution of the homogeneous problem (\ref{FBVP})
with fuzzy boundary conditions. $x_{cr}(t)$ can be computed by means of
analytical or numerical methods. Hence, (\ref{G_FBVP}) is reduced to solving
a homogeneous equation with fuzzy boundary conditions (\ref{FBVP}).
Therefore, we will investigate how to solve this problem.

We assume the solution $x_{un}$ of the problem (\ref{FBVP}) be a fuzzy set $%
\widetilde{X}$ of real functions such as $x(t)$. Each function $x(t)$ must
satisfy the differential equation and must have boundary values $a$ and $b$
from the sets $\widetilde{a}$ and $\widetilde{b}$, respectively. We define
the possibility (membership) of the function $x(t)$ to be equal to the least
possibility of its boundary values.

Mathematically, the fuzzy solution set can be defined as follows:\smallskip 

$\widetilde{X}=\{x(t)\mid x^{\prime \prime }+a_{1}(t)x^{\prime
}+a_{2}(t)x=0;\ $%
\begin{equation}
x(0)=a;\ x(T)=b;\ a\in \widetilde{a};\ b\in \widetilde{b}\}  \label{FS}
\end{equation}%
with membership function 
\begin{equation}
\mu _{\widetilde{X}}(x(t))=\min \left\{ \mu _{\widetilde{a}}(a),\mu _{%
\widetilde{b}}(b)\right\}
\end{equation}%
The solution $\widetilde{X}$, defined above, can be interpreted as a fuzzy
bunch of functions.

One can also interpret that we consider a FBVP as a set of crisp BVPs whose
boundary values belong to the fuzzy sets $\widetilde{a}$ and $\widetilde{b}$.

\subsection{A matrix representation of the solution in the crisp case}

Here we consider crisp BVP for second order homogeneous linear differential
equation:

\begin{equation}
\left\{ 
\begin{array}{c}
x^{\prime \prime }+a_{1}(t)x^{\prime }+a_{2}(t)x=0 \\ 
x(0)=a \\ 
x(T)=b%
\end{array}%
\right.  \label{CBVP}
\end{equation}

Let $x_{1}(t)$ and $x_{2}(t)$ be linear independent solutions of the
differential equation. Then the general solution is $%
x(t)=c_{1}x_{1}(t)+c_{2}x_{2}(t)$. For $c_{1}$ and $c_{2}$ we have the
following linear system 
\begin{equation}
\left\{ 
\begin{array}{c}
c_{1}x_{1}(0)+c_{2}x_{2}(0)=a \\ 
c_{1}x_{1}(T)+c_{2}x_{2}(T)=b%
\end{array}
\right.  \label{const}
\end{equation}

Below we obtain a matrix representation for the solution of the BVP. We
rewrite the linear system (\ref{const}) in matrix form:%
\begin{equation*}
M\ \mathbf{c}=\mathbf{u}
\end{equation*}%
where $M=\left[ 
\begin{array}{ll}
x_{1}(0) & x_{2}(0) \\ 
x_{1}(T) & x_{2}(T)%
\end{array}%
\right] $; $\mathbf{c}=\left[ 
\begin{array}{l}
c_{1} \\ 
c_{2}%
\end{array}%
\right] $; $\mathbf{u}=\left[ 
\begin{array}{l}
a \\ 
b%
\end{array}%
\right] $.

The solution of the linear system is 
\begin{equation}
\mathbf{c}=M^{-1}\ \mathbf{u}  \label{c}
\end{equation}%
We constitute a vector-function of linear independent solutions $\mathbf{s}%
(t)=\left[ x_{1}(t)\ \ \ x_{2}(t)\right] $. Then the general solution can be
rewritten in matrix form as%
\begin{equation*}
x(t)=\left[ x_{1}(t)\ \ \ x_{2}(t)\right] \left[ 
\begin{array}{l}
c_{1} \\ 
c_{2}%
\end{array}%
\right] =\mathbf{s}(t)\ \mathbf{c}
\end{equation*}

Using (\ref{c}) we have $x(t)=\mathbf{s}(t)\ M^{-1}\ \mathbf{u,}$ or, 
\begin{equation}
x(t)=\mathbf{w}(t)\ \mathbf{u=}w_{1}(t)\ a+w_{2}(t)\ b  \label{main}
\end{equation}%
where 
\begin{equation}
\mathbf{w}(t)=\mathbf{s}(t)\ M^{-1}  \label{w}
\end{equation}

\subsection{The solution method for FBVP}

Now we show how to find $\widetilde{X}(t)$ (the value of the solution for
the problem (\ref{FBVP}) at a time $t$).

Let linear independent solutions of the crisp equation (\ref{FBVP}),$\
x_{1}(t)$ and $x_{2}(t),$ be known. Then we can constitute the vector $%
\mathbf{w}$ (see, formula (\ref{w})). According (\ref{FS}) and (\ref{main})
we have: 
\begin{equation}
\widetilde{X}=\left\{ x(t)=\mathbf{w}(t)\ \mathbf{u}\mid \mathbf{u}=\left[ 
\begin{array}{l}
a \\ 
b%
\end{array}%
\right] ;\ a\in \widetilde{a};\ b\in \widetilde{b}\right\}   \label{FS_2}
\end{equation}

Consider a fixed time $t$. Put $\mathbf{v}=\mathbf{w}(t)$. Then from (\ref%
{FS_2}) we have: 
\begin{equation}
\widetilde{X}(t)=\left\{ \mathbf{v}\ \mathbf{u}\mid \mathbf{u}=\left[ a\ \ \
b\right] ^{T};\ a\in \widetilde{a};\ b\in \widetilde{b}\right\}
\end{equation}

To determine how is the set $\widetilde{X}(t)$ we consider the
transformation $T(\mathbf{u)=v\ u}$ (here $\mathbf{v}$ is a fixed vector).
One can see that $T:R^{2}\rightarrow R^{1}$ is a linear transformation.
Therefore, $\widetilde{X}(t)$ is the image of the set $\widetilde{B}=\left\{ 
\mathbf{u}=\left[ a\ \ \ b\right] ^{T}\mid a\in \widetilde{a};\ b\in 
\widetilde{b}\right\} =(\widetilde{a},\ \widetilde{b})$ under the linear
transformation $T(\mathbf{u)}$.

We remember some properties of linear transformations \cite{AR05}:

1. A linear transformation maps the origin (zero vector) to the origin (zero
vector).

2. Under a linear transformation the images of a pair of similar figures are
also similar.

3. Under a linear transformation the images of nested figures are also
nested.

In addition, we shall reference a property of fuzzy number vectors.

4. The fuzzy set $\widetilde{B}=(\widetilde{a},\ \widetilde{b})$ forms a
fuzzy region in the $ab$-coordinate plane, vertex of which is located at the
origin and boundary of which is a rectangle. Furthermore, the $\alpha $-cuts
of the region are rectangles nested within one another.

The facts 1-4 allow us to derive the following conclusion. The vector $%
\widetilde{B}$,\ components of which are the boundary values $\widetilde{a}$\
and $\widetilde{b}$, form a fuzzy rectangle in the $ab$-coordinate plane.
The linear transformation $T(\mathbf{u)}$ maps this fuzzy rectangle to a
fuzzy interval in $R^{1}$. Therefore, the solution at any time forms a fuzzy
number.

\subsection{Particular case when boundary values are triangular fuzzy numbers%
}

In particular, if $\widetilde{a}$\ and $\widetilde{b}$ triangular fuzzy
numbers, the $\alpha $-cuts of the region $\widetilde{B}=(\widetilde{a},\ 
\widetilde{b})$\ are nested rectangles, furthermore, they are similar.
According to the discussion above, their images are intervals that also are
nested and similar, consequently, form a triangular fuzzy number $\widetilde{%
X}(t)$. Therefore, $\widetilde{X}(t)$\ can be represented in the form $%
\widetilde{X}(t)=(\underline{x}(t),0,\overline{x}(t))$.\ Now we investigate
how to calculate $\underline{x}(t)$ and $\overline{x}(t)$.

Let $\widetilde{a}=(\underline{a},0,\overline{a})$,\ $\widetilde{b}=(%
\underline{b},0,\overline{b})$ and $\mathbf{w}(t)=(w_{1}(t),w_{2}(t))$.\
Since $\overline{x}(t)$ is the maximum value among all products $\mathbf{w}%
(t)\mathbf{\cdot u=}aw_{1}(t)+bw_{2}(t)$, we have:\smallskip 

$\overline{x}(t)=\max \left\{ \underline{a}w_{1}(t),\overline{a}%
w_{1}(t)\right\} +\max \left\{ \underline{b}w_{2}(t),\overline{b}%
w_{2}(t)\right\} $

$\underline{x}(t)=\min \left\{ \underline{a}w_{1}(t),\overline{a}%
w_{1}(t)\right\} +\min \left\{ \underline{b}w_{2}(t),\overline{b}%
w_{2}(t)\right\} $\smallskip 

Note that an $\alpha $-cut of $\widetilde{X}(t)$ can be determined by
similarity:%
\begin{equation*}
X_{\alpha }(t)=\left[ \underline{x_{\alpha }}(t),\overline{x_{\alpha }}(t)%
\right] =(1-\alpha )\left[ \underline{x}(t),\overline{x}(t)\right]
=(1-\alpha )X_{0}(t)
\end{equation*}

Formulas for $\underline{x}(t)$ and $\overline{x}(t)$ allow us to represent
the solution in a new way: 
\begin{equation*}
\widetilde{X}(t)=w_{1}(t)\ \widetilde{a}+w_{2}(t)\ \widetilde{b}
\end{equation*}%
where the operations assumed to be multiplication of real number with fuzzy
one, and addition of fuzzy numbers.

\subsection{General case when boundary values are parametric fuzzy numbers}

In the general case, when $\widetilde{a}$\ and $\widetilde{b}$ are arbitrary
fuzzy numbers, the solution can be obtained by using $\alpha $-cuts. Let $%
a_{\alpha }=\left[ \underline{a_{\alpha }},\overline{a_{\alpha }}\right] $
and $b_{\alpha }=\left[ \underline{b_{\alpha }},\overline{b_{\alpha }}\right]
$. Then $B_{\alpha }=\left[ \underline{a_{\alpha }},\overline{a_{\alpha }}%
\right] \times \left[ \underline{b_{\alpha }},\overline{b_{\alpha }}\right] $%
. By similar argumentation to the preceding case, for the $\alpha $-cut of
the solution we obtain the following formulas:$\smallskip $

$X_{\alpha }(t)=\left[ \underline{x_{\alpha }}(t),\overline{x_{\alpha }}(t)%
\right] \smallskip $

\noindent where$\smallskip $

\noindent $\overline{x_{\alpha }}(t)=\max \left\{ \underline{a_{\alpha }}%
w_{1}(t),\overline{a_{\alpha }}w_{1}(t)\right\} +\max \left\{ \underline{%
b_{\alpha }}w_{2}(t),\overline{b_{\alpha }}w_{2}(t)\right\} $

\noindent $\underline{x_{\alpha }}(t)=\min \left\{ \underline{a_{\alpha }}%
w_{1}(t),\overline{a_{\alpha }}w_{1}(t)\right\} +\min \left\{ \underline{%
b_{\alpha }}w_{2}(t),\overline{b_{\alpha }}w_{2}(t)\right\} \smallskip $

Based on the formulas above we can conclude, that the representation for
solution 
\begin{equation}
\widetilde{X}(t)=w_{1}(t)\ \widetilde{a}+w_{2}(t)\ \widetilde{b}  \label{mf}
\end{equation}%
is valid in general. Thus, the solution, defined by formula (\ref{FS}),
becomes the same as the solution obtained from (\ref{main}) by extension
principle.

\textit{Remark}: The approach is valid also for the general case, when $n$th-order 
$m$-point boundary value problem is considered.

\subsection{Solution algorithm}

Based on the arguments above, we propose the following algorithm to solve
the problem (\ref{G_FBVP}):

1. Represent the boundary values as $\widetilde{A}=a_{cr}+\widetilde{a}$\
and $\widetilde{B}=b_{cr}+\widetilde{b}$.

2. Find linear independent solutions $x_{1}(t)$ and $x_{2}(t)$ of the crisp
differential equation $x^{\prime \prime }+a_{1}(t)x^{\prime }+a_{2}(t)x=0$.
Constitute the vector-function $\mathbf{s}(t)=\left[ x_{1}(t)\ \ \ x_{2}(t)%
\right] $, the matrix $M$ and calculate the vector $\mathbf{w}%
(t)=(w_{1}(t),w_{2}(t))$ by formula (\ref{w}).

3. Find the solution $x_{cr}(t)$ of the non-homogeneous crisp problem (\ref%
{N-H_CBVP}).

4. The solution of the given problem (\ref{G_FBVP}) is 
\begin{equation}
\widetilde{x}(t)=x_{cr}(t)+w_{1}(t)\ \widetilde{a}+w_{2}(t)\ \widetilde{b}
\end{equation}

\section{Examples}

\noindent\textbf{Example 1}. Solve the FBVP: 
\begin{equation}
\left\{ 
\begin{array}{c}
x^{\prime \prime }-3x^{\prime }+2x=4t-6 \\ 
x(0)=(1.5,\ 2,\ 3) \\ 
x(1)=(\ \ 2,\ \ 3,\ 4)%
\end{array}%
\right.  \label{Ex1}
\end{equation}

\noindent\textbf{Solution}: We represent the solution as $\widetilde{x}%
(t)=x_{cr}(t)+\widetilde{x}_{un}(t)$.

1. We solve crisp non-homogeneous problem

$\left\{ 
\begin{array}{c}
x^{\prime \prime }-3x^{\prime }+2x=4t-6 \\ 
x(0)=2 \\ 
x(1)=3%
\end{array}%
\right. $\newline
and find the crisp solution \newline
$x_{cr}(t)=2t+\frac{1}{e^{2}-e}\left[ 2(e^{2+t}-e^{1+2t})+(e^{2t}-e^{t})%
\right] $ (the thick line in Fig.1). 
\begin{figure}[h]
\centerline{\epsfig{figure=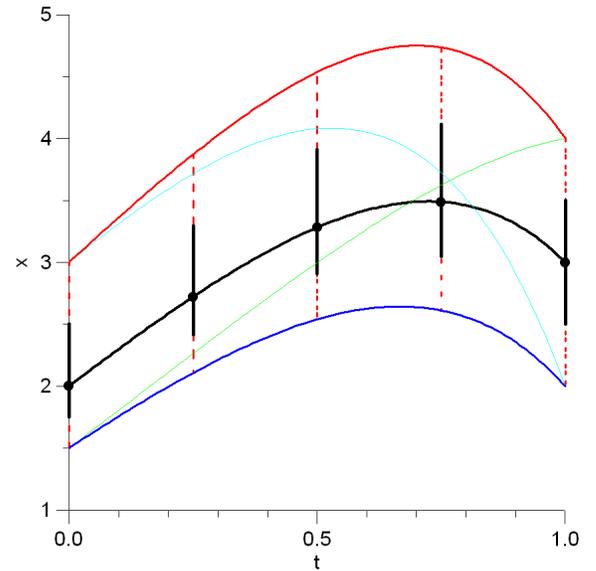,width=75 mm, angle=0} }
\caption{The fuzzy solution and its $\protect\alpha =0.5$-cut  for Example 1. Dashed and thick bars represent the values of the fuzzy solution and its $\protect\alpha =0.5$-cut at different
times, respectively.}
\label{fig:fig1}
\end{figure}

2. We consider fuzzy homogeneous problem

$\left\{ 
\begin{array}{c}
x^{\prime \prime }-3x^{\prime }+2x=0 \\ 
x(0)=(-0.5,\ 0,\ 1) \\ 
x(1)=(-1,\ 0,\ 1)%
\end{array}
\right. $

$x_{1}(t)=e^{t}$ and $x_{2}(t)=e^{2t}$ are linear independent solutions for
the differential equation. Then $\mathbf{s}(t)=\left[ e^{t}\ \ \ e^{2t}%
\right] $, \newline $M=\left[ 
\begin{array}{ll}
1 & 1 \\ 
e & e^{2}%
\end{array}%
\right] $ \ and \newline
$\mathbf{w=s}(t)\ M^{-1}=\frac{1}{e^{2}-e}\left[ e^{2+t}-e^{1+2t}\ \ \
e^{2t}-e^{t}\right] $.\smallskip 

\noindent The formula (\ref{mf}) gives the solution of homogeneous problem:
\smallskip 

$\widetilde{x}_{un}(t)=\frac{1}{e^{2}-e}((e^{2+t}-e^{1+2t})\ (-0.5,\ 0,\ 1)+$%
\begin{equation}
(e^{2t}-e^{t})\ (-1,\ 0,\ 1))
\end{equation}%
where the arithmetic operations are considered to be fuzzy operations. We
add this solution to the crisp solution and get the fuzzy solution of the
given FBVP (\ref{Ex1}):\smallskip 

$\widetilde{x}(t)=2t+\frac{1}{e^{2}-e}((e^{2+t}-e^{1+2t})\ (1.5,\ 2,\ 3)+$%
\begin{equation}
(e^{2t}-e^{t})\ (0,\ 1,\ 2))
\end{equation}

The fuzzy solution $\widetilde{x}(t)$ forms a band in the $tx$-coordinate
plane (Fig. 1). Since $w_{1}(t)>0$ and $w_{2}(t)>0$ for $0<t<T$, the upper
border of the band, $\overline{x}(t)$, becomes the solution of the crisp
non-homogeneous problem with the upper boundary values $\overline{A}=3$ and $%
\overline{B}=4$, while the lower border $\underline{x}(t)$ corresponds to $%
\underline{A}=1.5$ and $\underline{B}=2$: 
\begin{eqnarray*}
\overline{x}(t) &=&2t+\frac{1}{e^{2}-e}((e^{2+t}-e^{1+2t})\cdot
\ 3\ +(e^{2t}-e^{t})\cdot 2) \\
\underline{x}(t) &=&2t+\frac{1}{e^{2}-e}((e^{2+t}-e^{1+2t})\cdot
1.5+(e^{2t}-e^{t})\cdot 0)
\end{eqnarray*}

We can express the solution $\widetilde{x}(t)$ also via $\alpha $-cuts,
which are intervals $x_{\alpha }(t)=\left[ \underline{x_{_{\alpha }}}(t),\ 
\overline{x_{_{\alpha }}}(t)\right] $ at any time $t$. Since the boundary
values are triangular fuzzy numbers, $\widetilde{x}_{un}(t)$ also is a
triangular fuzzy number, say $\widetilde{x}_{un}(t)=(\underline{x_{un}}(t),\
0,\ \overline{x_{un}}(t))$. Consequently, an $\alpha $-cut of $\widetilde{x}%
_{un}(t)$ can be determined by similarity with coefficient $(1-\alpha )$,
i.e.%
\begin{equation*}
x_{un,\ \alpha }(t)=(1-\alpha )\left[ \underline{x_{un}}(t),\ \overline{%
x_{un}}(t)\right] 
\end{equation*}%
Adding the crisp solution gives the $\alpha $-cut of the solution $%
\widetilde{x}(t)$:\smallskip 

\noindent $\left[ \underline{x_{_{\alpha }}}(t),\ \overline{x_{_{\alpha }}}%
(t)\right] =2t+(1-\alpha )\cdot \frac{1}{e^{2}-e}((e^{2+t}-e^{1+2t})\ \left[
1.5,\ 3\right] +$%
\begin{equation*}
(e^{2t}-e^{t})\ \left[ 0,\ 2\right] )
\end{equation*}

In Fig. 1 we show the fuzzy solution (dashed bars) and its $%
\alpha =0.5$-cut (thick bars) at different times.\smallskip 

\noindent\textbf{Example 2}. Solve the FBVP: 
\begin{equation}
\left\{ 
\begin{array}{c}
x^{\prime \prime }+16x=47-8t^{2} \\ 
x(0)=(\ \ 2,\ \ 3,\ 3.5) \\ 
x(2)=(0.5,\ 1,\ 1.5)%
\end{array}
\right. \label{Ex2}
\end{equation}

\noindent\textbf{Solution}:

Associated crisp non-homogeneous problem

$\left\{ 
\begin{array}{c}
x^{\prime \prime }+16x=47-8t^{2} \\ 
x(0)=3 \\ 
x(2)=1%
\end{array}
\right. $

\noindent has the solution $x_{cr}(t)=3-0.5t^{2}$ (thick line in Fig. 2).

To find the uncertain part of the fuzzy solution, $\widetilde{x}_{un}(t)$,
we solve fuzzy homogeneous problem

$\left\{ 
\begin{array}{c}
x^{\prime \prime }+16x=0 \\ 
x(0)=(\ -1,\ 0,\ 0.5) \\ 
x(2)=(-0.5,\ 0,\ 0.5)%
\end{array}
\right. $

\begin{figure}[h]
\centerline{\epsfig{figure=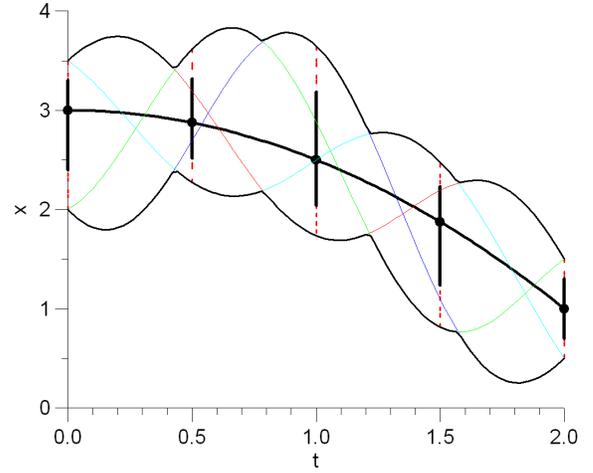,width=75 mm, angle=0} }
\caption{The fuzzy solution and its $\protect\alpha =0.6$-cut for Example 2. Dashed and thick bars represent the values of the fuzzy solution and its $\protect\alpha =0.6$-cut at different
times, respectively.}
\label{fig:fig2}
\end{figure}

$x_{1}(t)=\cos 4t$ and $x_{2}(t)=\sin 4t$ are linear independent solutions
for the differential equation. Then

\noindent $\mathbf{s}(t)=\left[ \cos 4t\ \ \ \sin 4t\right] $, \ $M=\left[ 
\begin{array}{ll}
1 & 0 \\ 
\cos 8 & \sin 8%
\end{array}%
\right] $ \ and \ \newline
$\mathbf{w=s}(t)\ M^{-1}=\frac{1}{\sin 8}\left[ \sin (8-4t)\ \ \ \sin 4t%
\right] $.

Using the formula (\ref{mf}) we obtain the solution of homogeneous problem
and adding the crisp solution we get the solution of the given FBVP (\ref%
{Ex2}): \smallskip 

$\widetilde{x}(t)=3-0.5t^{2}+\frac{1}{\sin 8}(\sin (8-4t)\ (-1,\ 0,\ 0.5)+$

\begin{equation}
\sin 4t\ (-0.5,\ 0,\ 0.5))  \label{SolEx2}
\end{equation}

Fuzzy solution generates a band in $tx$-plane (Fig. 2). Unlike Example 1,
the functions $w_{1}(t)$ and $w_{2}(t)$ takes both
positive and negative values in the interval $0<t<T$. Because of that, in generation of upper and lower
borders of the band, $\overline{a}$ and $\underline{a}$, $\underline{b}
$ and $\overline{b}$ take charge in alternately.

\bigskip

\section{Conclusion}

In this paper we have investigated the fuzzy boundary value problem as a set of crisp
problems. We have proposed a solution method based on the properties of linear transformations. For clarity we have explained the
proposed method for second order linear differential equation. We have shown that the fuzzy solution by our method coincides with the
solution by extension principle. We are planning to make a comparative analysis between
the proposed method and the method including generalized Hukuhara
derivative in future.




\begin{thebibliography}{99}
\bibitem{Bede06} B. Bede, "A note on "two-point boundary value problems
associated with non-linear fuzzy differential equations"", Fuzzy Sets and
Systems, 157 (2006) 986-989.

\bibitem{BG05} B. Bede and S.G. Gal, "Generalizations of the
differentiability of fuzzy number valued functions with applications to
fuzzy differential equation", Fuzzy Sets and Systems, 151 (2005) 581-599.

\bibitem{CCRF06} Y. Chalco-Cano and H. Román-Flores, "On the new solution of
fuzzy differential equations", Chaos, Solitons, Fractals, 38 (2006) 112-119.

\bibitem{BRB07} B. Bede, I.J. Rudas, and A.L. Bencsik, "First order linear
fuzzy differential equations under generalized differentiability", Inform.
Sci., 177 (2007) 1648-1662.

\bibitem{CCRFRM08} Y. Chalco-Cano, H. Román-Flores, and M.A. Rojas-Medar,
"Fuzzy differential equations with generalized derivative", in: Proc. 27th
NAFIPS Internat. Conf. IEEE, 2008.

\bibitem{CCRMRF07} Y. Chalco-Cano, M.A. Rojas-Medar, and H. Román-Flores,
"Sobre ecuaciones diferenciales difusas", Bol. Soc. Española Mat. Aplicada,
41 (2007) 91-99.

\bibitem{KhN10} A. Khastan and J.J. Nieto, "A boundary value problem for
second order fuzzy differential equations", Nonlinear Anal., 72 (2010)
3583-3593.

\bibitem{BF00} J.J. Buckley and T. Feuring, "Fuzzy differential equations",
Fuzzy Sets and Systems, 110 (2000) 43-54.

\bibitem{BF01} J.J. Buckley and T. Feuring, "Fuzzy initial value problem for
Nth-order linear differential equation", Fuzzy Sets and Systems, 121 (2001)
247-255.

\bibitem{Hul97} E. Hüllermeier, "An approach to modeling and simulation of
uncertain dynamical systems", Internat. J. Uncertainty, Fuzziness,
Knowledge-Based Systems, 5 (1997) 117-137.

\bibitem{MBCRB07} M. Misukoshi, L.C. Barros, Y. Chalco-Cano, H. Romá%
n-Flores, and R.C. Bassanezi, "Fuzzy differential equations and the
extension principle", Inform. Sci., 177 (2007) 3627-3635.

\bibitem{GAF11} N.A. Gasilov, Ş.E. Amrahov, and A.G. Fatullayev, "A
geometric approach to solve fuzzy linear systems of differential equations",
Appl. Math. Inf. Sci., 5 (2011) 484-495.

\bibitem{AR05} H. Anton and C. Rorres, Elementary Linear Algebra,
Applications Version: 9th Edition, John Wiley \& Sons, 2005.
\end{thebibliography}
\end{document}